\documentclass[10pt]{amsart}
\usepackage{amsfonts,amssymb,amscd}
\usepackage{pstricks}
\usepackage{pst-plot}
\usepackage{mathrsfs}

\let\mathcal\mathscr

\usepackage[all,ps,cmtip,color]{xy}

\newtheorem{theo}{Theorem}
\newtheorem{coro}{Corollary}
\newtheorem{prop}{Proposition}
\newtheorem{lemm}{Lemma}

\begin{document}

\def\om{\omega}
\def\ot{\otimes}
\def\we{\wedge}
\def\wec{\wedge\cdots\wedge}
\def\op{\oplus}
\def\ra{\rightarrow}
\def\lra{\longrightarrow}
\def\fso{\mathfrak so}\def\fg{\mathfrak g}
\def\cO{\mathcal{O}}
\def\cS{\mathcal{S}}
\def\fsl{\mathfrak sl}
\def\fb{\mathfrak b}
\def\fh{\mathfrak h}
\def\fn{\mathfrak n}
\def\PP{\mathbb P}\def\QQ{\mathbb Q}\def\ZZ{\mathbb Z}\def\CC{\mathbb C}
\def\RR{\mathbb R}\def\HH{\mathbb H}\def\OO{\mathbb O}
\def\smc{\cdots }
\def\Gr{{\mathcal Gr}}
\def\Gm{\mathbf{G_m}}
\def\tP{{\tilde{P}}}
\newcommand{\comb}[2]{\binom{#1}{#2}}

\title{Quantum cohomology and  the Satake isomorphism}

\author{V. Golyshev, L. Manivel}

\begin{abstract}

We prove that the geometric Satake correspondence admits quantum corrections for minuscule
Grassmannians of Dynkin types $A$ and $D$. We find, as a corollary, that the quantum connection of a
spinor variety $OG(n,2n)$
can be obtained as the half-spinorial representation of that of the quadric $Q_{2n-2}$.
We view the (quantum)
cohomology of these Grassmannians as  endowed
simultaneously with two structures, one of a module over the algebra of symmetric functions,
and the other, of a module over the Langlands dual Lie algebra,
and investigate the interaction between the two. In particular, we study primitive classes $y$
in the cohomology of a minuscule Grassmannian $G/P$ that are
characterized by the condition that the operator of cup product by $y$
is in the image of the Lie algebra action.
Our main result states that quantum correction preserves primitivity.
We provide a quantum counterpart to a
result obtained by V. Ginzburg in the classical setting by giving explicit
formulas
for the quantum corrections to homogeneous primitive elements.

\end{abstract}

\maketitle

\section{Motivations}

The statements of this paper can be viewed as suggested by the classical Tannakian philosophy applied in the quantum setup.
We refer the reader to an  excellent textbook \cite{and} which explains in detail how the motivic approach works in
the classical setup by turning, or at least trying to turn, every object into a representation of a certain Tannakian group.
The classical realisations and the classical
motivic group are meant to treat essentially transcendental pieces of algebraic varieties
but would be too small if we were to treat varieties, say, with cellular decomposition.
Too many correspondences would imply too many projectors, and the variety would split up into a direct
sum of Tate motives. Is it possible, for varieties such as cellular ones, to
translate the incidence structure of cells, or cycles, adequately into a structure of a module over a group?
Is it possible to do
that effectively and explicitly, approaching the problem of classifying the incidences of algebraic geometry
by that of classifying certain representations?

Where a linear
calculus of incidence is sought for, quantum motives are a seemingly feasible substitute.

\subsection{Classical semisimple motives and the motivic group.}
In his 2008 IHES talk, Manin \cite{man-ihes} views the structure of an algebra over the cyclic modular operad
$\coprod_g h(\overline{M}_{g,n})$ on a total motive
as  the motivic core of quantum cohomology, thus emphasizing the necessity to work with all genera.
If one accepts the idea of downshifting to small genus $0$
quantum motives, one finds today's situation strikingly similar, in many respects, to the one
with the classical motives  in the early 70'ies.
Consider algebraic varieties over $\QQ$ and a Weil cohomology theory, such as the Betti cohomology
of the underlying analytic spaces as $\ZZ$--lattices. How can one detect algebraic morphisms at the level of
the maps between the cohomology groups?
The ultimate hope was that it might be possible
to turn $H(V)$'s into a module over a `reductive group' $M$
such that

\begin{enumerate}

\item any algebraic morphism  $h: H(V_1) \longrightarrow H(V_2)$ is $M$--invariant;
and any $M$--invariant algebraic morphism  $h: H(V_1) \longrightarrow H(V_2)$ is algebraic;

\item any irreducible $M$--representation is contained in one that comes from geometry,

\end{enumerate}
thus being able to classify the essentially transcendental pieces of algebraic varieties
as representations of a
certain group.

Different versions of this formulation are possible: one could imagine an algebraic group $M_V$ that controls algebraicity specifically in $V$ and all of its cartesian powers; or a pro--algebraic $M$ that would control
 all smooth varieties $V$ simultaneously.

In either case, having required as much as above, one has essentially forced oneself
into the scheme of
classical semisimple (homological) motives. \nocite{man-corr}
Each $M$--module  decomposes into a direct sum of submodules and each projector is algebraic. Hence, the construction must start with adjoining the category of algebraic varieties (or rather its opposite)
by the images of the projectors in the additive category of correspondences.

As is well known, this dream never fully came true, the obstacle being the standard motivic conjectures.
If, for example,
we take for $H$ the $l$--adic cohomology of the variety over $\overline{\QQ}$ with its Galois action, (1) above
is essentially the Tate conjecture, and (2) is the Fontaine--Mazur conjecture.

\subsection{More correspondences.} A way out suggested by Yves Andr\'e (see \cite{and}) is
to substitute the Hodge group for the  motivic group
thereby passing to Hodge motives.
Recall that the Hodge group of a (smooth complex projective) variety $V$
is defined to be the minimal algebraic $\mathbb{Q}$--group $Hdg_V$ such that $Hdg_V(\mathbb{R})$
contains the image of the Hodge
action $h: S^1 \longrightarrow GL(H(V) \otimes \mathbb{R})$. The Hodge classes are the classes
invariant under the action of the Hodge group.
Consider the category obtained by adjoining the images of projectors with respect to the Hodge correspondences (rather than algebraic) to obtain Hodge motives.
It is clear that there is a Hodge realisation functor from Hodge motives to Hodge structures.

More generally, with any realization  in a [semisimple]
Tannakian category $\mathcal{T}$ comes a similar construction where the images of projectors that are invariant
with respect to the action of the Tannakian group $T_\mathcal{T}$ are adjoined. However, we want to apply
this idea in a situation where the Tannakian group is as big as to mix up
cohomologies in different dimensions.

\subsection{Fewer correspondences.}
One is led, therefore, to a motivic construction with fewer projectors, namely those that are
`neutral with respect to the Gromov--Witten calculus', a notion exemplified in \cite{BC-FK}.
A minimal possible construction is probably \emph{homological quantum motives} in the spirit of Andr\'e, with respect to the D--module realization. Let $F$ be a Fano variety;
we will assume for a while that
$\mathop{\mathrm{Pic}} F = \ZZ$, and that all classes $HF$ are algebraic. Recall that
 a three--point correlator
$\langle a, b ,c \rangle _d$ of three effective algebraic classes in $H(F)$
is the number of rational curves of $-K_F$--degree $d,$
intersecting the generic representatives of these classes. The algebra $QH(F)$
is $\CC[q,q^{-1}]$--module isomorphic to $HF
\otimes \CC [q,q^{-1}]$ with the following multiplication structure.
Denote the classes by  $A, B, C$
so that $A = a \otimes 1$, etc.  Put
$$A \cdot B = \sum_{c , d \ge 0}
\langle a,b,c ^\vee \rangle _d
\; q^d C , $$ where $C ^\vee$-- is the class dual to
$C$ and the summation w.r. to $c$
is over some basis of $HF$.
 Interpret $QHF$ as the space of sections of the constant vector bundle over
$\Gm = \mathop{\mathrm{Spec}} \CC[q,q^{-1}]$.
Put $D = q\frac{d}{dq}$ and denote by
$M$ the operator of quantum multiplication by $-K_F$ expressed in the basis $\xi$ of horizontal
sections (i.e. those of the shape $a \otimes 1$).

The quantum connection
$D \xi = \xi M$ turns the cohomology of $F$ into an object of the Tannakian category
$DE(\Gm)$. Recall that the differential fundamental group of a variety $V$ is defined to be
the group of the automorphisms of the fiber functor of the
Tannakian category of differential equations on $V$.
In this parlance, the vector space $HF$ acquires a structure of a
$\pi_1^\text{diff}(\Gm)$--module.

\subsection{Group actions and the effectivity problem.} Quantum multiplication is respected by automorphisms:
for any $g \in G$ acting on a Fano $F$ we have $g^* a \cdot g^*b = g^*(a \cdot b)$. In particular,
the Schur quantum projectors arise in the cohomology of the Cartesian powers of $F$.
Though defined as a quantum
motive, the image of the respective projector $\pi$
may in principle be effective at the level of the realization, that is, there may exist a Fano or a Fano--like
object $F^\pi$ such that $S^\pi QHF)= QHF^\pi$ as DE's.


\subsection{Problem.}  Tell when a given quantum motive is effective.

\newcommand{\p}{\mathbf{p}}
\renewcommand{\P}{{\mathbb{P}}}
\newcommand{\G}{{\mathbb{G}}}

\subsection{Grassmannians: Schubert and Satake.} \label{Schubert-and-Satake}
This is an instance of effectivity.
Denote $\P = \P^{n-1}=G(1,n)$ and $\G = G(l,n)$. Then:

\begin{enumerate}

\item the quantum connection of $\P$ is irreducible;

\item the  $l$-th wedge power of the quantum connection of $\P$ is the quantum connection of $\G$.

\end{enumerate}

This theorem may be viewed as an instance of a more general `quantum Satake' principle that says, roughly,
that the geometric Satake correspondence of Beilinson--Drinfeld--Lusztig--Ginzburg--Mirkovic--Vilonen
\cite{MV07}, \cite{BF08}
admits quantum corrections for minuscule Grassmannians. Before
we review the geometric Satake in the next section, let us take a look at how the Schubert and the Satake
structures interplay in the cohomology of Grassmanians of Dynkin type A. A theorem of V.~Ginzburg
\cite[1.3.2]{Ginzburg08} provides a description in the general case (and it is essentially
its quantum counterpart that we study in this paper).

\begin{itemize}

\item \bf Schubert. \rm By definition, the \emph{Schubert structure} on the cohomology of $\G$ is the map
$$\mathrm{Sch}_\G: \Lambda \longrightarrow H\G$$
from the ring of symmetric functions.  The kernel ideal is
$$I=\langle e_{l+1},e_{l+2}, \dots h_{n-l+1}, h_{n-l+1},\dots, h_n\rangle .$$
Similarly, the cohomology of $\P$ is endowed with the Schubert structure $\mathrm{Sch}_\P: \Lambda \longrightarrow H\P$.

\medskip

\item \bf Satake.  \rm The \emph{Satake structure}
$\mathrm{Sat}: \Lambda^lH\P\ \longrightarrow H\G$
is the identification of the cohomology of $\G$ as a minuscule Schubert cell in the affine Grassmannian
of $GL(n)$
with the wedge power of cohomology of $\P$, also interpreted as a minuscule Schubert cell
in the same affine Grassmannian. Concretely, the Satake identification
is realized by the `take the span' correspondence between $\P ^l$ and $\G$; one has
$$\mathrm{Sat}:
\sigma_{\lambda_1+l-1} \wedge \sigma_{\lambda_2+l-2}\wedge \dots \wedge \sigma_{\lambda_l}
\mapsto
\sigma_\lambda$$

The geometric Satake principle of Ginzburg  predicts, in particular, that $\mathrm{Sat}$ extends to endow $H\G$ with a natural  $gl(H\P)$--module structure,
and that the operator $\p_1 \in gl(H\P)$ of multiplication by $c_1(\P)$ acts in $H\G$
by multiplication by $c_1(\G)$. In our situation, this is easy to see directly.

\end{itemize}

\subsection{Quantum corrections to Schubert and Satake.} The structures above admit quantum deformations as follows:
\begin{itemize}
\item \bf Quantum Schubert. \rm By theorem of Siebert and Tian,
the Schubert structure on the quantum cohomology ring of $\G$ is again a map
$$\mathrm{Sch}_\G^q: \Lambda [q] \longrightarrow QH\G$$
The kernel ideal is now
$$I_q= \langle e_{l+1},e_{l+2}, \dots h_{n-l+1}, h_{n-l+1},\dots, h_n+(-1)^l q\rangle .$$
Similarly, the quantum cohomology ring of $\P$ is endowed with the Schubert structure
$\mathrm{Sch}_\P^q$.

\medskip

\item \bf Quantum Satake. \rm
Let now $\p_i^q \in gl(H\P,\CC [q])$ be the operator of quantum multiplication
by the class $p_i \in QH\P$. Then, upon the Satake identification, $\p_i^q$ acts in $QH\G$ by quantum multiplication
by the class $\mathrm{Sch}_\G(p_i)$ in $QH\G$.

\end{itemize}

\medskip

\subsection{Review of Ginzburg's theorem.} \label{review-Ginzburg}
The Satake isomorphism \cite{gi}
involves an affine Schubert variety $X_\lambda$ in the affine Grassmannian
$\Gr$ of a complex Lie group $G$, indexed by a one parameter subgroup $\lambda\in X_*(T)$
of a maximal torus $T$ of $G$, up to the action of the Weyl group. If $G^L$ is the Langlands dual
of $G$, with maximal torus $T^L$, then $X_*(T)=X^*(T^L)$ and $\lambda$ defines an irreducible
$G^L$-module $V^L_\lambda$. The claim is that $V^L_\lambda$ can be identified with the intersection
cohomology of $X_\lambda$, in such a way that the natural action of $H^*(\Gr,\CC)$ on $IH(X_\lambda)$
is identified with the action of $U(\mathfrak{a})$ on $V^L_\lambda$, where $\mathfrak{a}\subset\fg^L$
is the centralizer of a regular nilpotent element in the Lie algebra $\fg^L$ of $G^L$.
Once a Cartan subalgebra $\fh$ has been chosen, one can for example chose the regular nilpotent element
$$x=X_{-\alpha_1}+\cdots +  X_{-\alpha_r},$$
where $\alpha_1, \ldots , \alpha_r$ are the simple roots (with respect to some choice of
positive roots) and $X_{\alpha_i}$ belongs to the root space $\fg^L_{\alpha_i}$ of $\fg^L$.

\smallskip
In the case where $G$ is adjoint, $X_*(T)=Q^\vee$ is dual to the root lattice $Q$ and can be
identified with the weight lattice $P^L$ of $G^L$. The connected components
of $\Gr$ are indexed by the element of $P/Q$, the quotient of the weight lattice by the root lattice. Recall that the
nonzero elements in $P/Q\simeq Q^\vee/P^\vee\simeq P^L/Q^L$ are in natural correspondence with
the minuscule fundamental weights of $G^L$, which are the minimal dominant elements in the $Q$-classes
of $P$. In particular, if $\omega$ is such a weight, then $X_\omega=G/P_\omega$ is a projective $G$-orbit,
in particular it is smooth and the Satake isomorphism is an identification
$$H^*(G/P_\omega,\CC)\simeq V^L_\omega.$$

\medskip

\subsection{Problem.} Does the geometric Satake principle
admit quantum corrections for non--minuscule Schubert cells in $\Gr (GL(N))$? In other words,
is it true that $QIC(X_\lambda)$ is isomorphic to $S^\lambda QH(\P)$ as differential equations?

\smallskip

We remark that even a universally established definition of quantum cohomology of such spaces as
$X_\lambda$'s is lacking. In fact we believe that such a definition should come with an in--built
reference to its behavior within the totality of all $X_\lambda$'s.

\section{Statement of the main theorem}

Our main result describes quantum correction to the multiplication in the cohomology of a minuscule
$G/P_\omega$ by generalizing the interplay between the Schubert and the Satake structures
described in the classical case in \ref{Schubert-and-Satake}.

\subsection{Schubert's Cartan subalgebra.}
 For $s \in W/W_P$ the
respective Schubert class is the cohomology class of $\overline{BsP/P}$, where $P=P_\omega$.
The Satake structure is the structure of $\mathfrak{g}^L$--module on $H^*(G/P_\omega,\CC)$. The
elements of $\mathfrak{g}^L$ that diagonalize in the Schubert basis form a Cartan subalgebra which we
denote by $\mathfrak{h}$. The subalgebra $\mathfrak{h}$ comes with a choice of a positive chamber if we demand
that the regular nilpotent
$x=X_{-\alpha_1}+\cdots +  X_{-\alpha_r}$ act by multiplication by the hyperplane class of $G/P_\omega$.
Thus, the fundamental class becomes the highest weight vector.
Recall that we have denoted by $\mathfrak{a}$ the centralizer of $x$.

\subsection{Kostant's Cartan subalgebra.}

Define the {\it cyclic element}:
$$x_q=X_{-\alpha_1}+\cdots +  X_{-\alpha_r}+qX_{\psi},$$
where $\psi$ denotes the highest root. Then, by Kostant \cite{kos},
the cyclic element $x_q$ is regular and semisimple. In particular its centralizer
$\mathfrak{a}_q$ is a Cartan subalgebra of $\mathfrak{g}^L$.
The two Cartan subalgebras $\mathfrak{h}$ and $\mathfrak{a}_q$
are {\it in apposition}, according to the terminology of \cite{kos}.
There is an isomorphism
$$\pi_q : \mathfrak{a}\longrightarrow\mathfrak{a}_q,$$
defined through the decomposition $\mathfrak{g}=\mathfrak{n}^-
\oplus\mathfrak{h}\oplus\mathfrak{n}$, where $\mathfrak{n}$ (resp. $\mathfrak{n}^-$)
is the span of the positive (resp. negative) root spaces. Indeed, $\mathfrak{a}$
is contained in $\mathfrak{n}^-$, and the projection of $\mathfrak{a}_q$ to
 $\mathfrak{n}^-$ is an isomorphism onto $\mathfrak{a}$. Our $\pi_q$ is the
inverse isomorphism.

\subsection{} Define now the map $\theta$  as the composition
$$
\begin{CD}
\mathfrak{a}   @>\mathrm{Sat}>>   \mathop{\mathrm{End}} (H^*(G/P_\omega,\CC))
@>\mathrm{ev_1}>> H^*(G/P_\omega,\CC),
\end{CD}
$$
where $\mathrm{ev_1}$ is the evaluation morphism on the fundamental class.
The map $\theta$ is injective. We denote by $\theta_q$ the composition of
$\theta$ with the injection of $H^*(G/P_\omega,\CC)$ into $QH^*(G/P_\omega,\CC)$.

\begin{theo}
Let $y\in \mathfrak{a}$. Then the quantum product in $QH^*(G/P_\omega)$ by the
class $\theta_q(y)$ coincides with the action on  $V^L_\omega$ by the
element $\pi_q(y)$. In other words, the following diagram is commutative:
$$
\xymatrix{
\mathfrak{a} \ar[r]^{\pi_q}  \ar[d]^{\theta_q} &\mathfrak{a}_q \ar[d]_{\mathrm{Sat}}  \\
QH^*(G/P_\omega,\CC) \ar[r]^{\star\hspace*{6mm} }  & \mathop{\mathrm{End}} (QH^*(G/P_\omega,\CC)).
}
$$

In particular, the quantum product by the hyperplane class $h=\theta(x)$
coincides with the Lie action of $\pi_q(x)=x_q$.
\end{theo}

Thus, the second part of the statement is the minuscule version
of the quantum Chevalley formula.

\bigskip

In order to restate Theorem 1 more
explicitly, we take a closer look at the
structure of $\pi_q$.

\subsection{The projection $\pi_q$ and the extended power classes.}
For any exponent
$d$ of $\fg$, there exists an element $y_d\in\mathfrak{a}$ of the form
$$y_d=\sum_{ht(\alpha)=d}y_d^\alpha X_{-\alpha}.$$
This follows from Kostant's results and the interpretation of the sequence
of exponents as the dual partition to the partition defined by the number
of roots of different heights.
Correspondingly,
there exists a unique element $y_{d,q}=\pi_q (y_d)\in\mathfrak{a_q}$, whose projection to
$\fn^-$ is $y_d$. By homogeneity, it can be written in the form $y_{d,q}=y_d+qz_d$,
where $z_d\in\mathfrak{n}$ is a combination of root vectors of height $h-d$.

Denote by $n_i(\beta)$ the
coefficient of a positive root $\beta$ on the simple root $\alpha_i$.
Put $q(\beta)=\prod_in_i(\psi)^{n_i(\beta)}$.

By Kostant's duality theorem \cite[Theorem 6.7]{kos}, we have
$$z_d=\zeta_d\sum_{ht(\beta)=h-d}y_{h-d}^{\beta}q(\beta) X_{\alpha}.$$
Applying $\theta$ to $y_d$ defines a class $p_d=\theta(y_d)$ in $H^{2d}(G/P_\omega,\CC)$, which
equals the result of applying $y_d$ to the fundamental class.
We will see shortly that the classes $p_d$ should be interpreted
as the analogs of the power sum classes in the cohomology of classical grassmannians.
We can define inductively a weight basis of $V^L_\omega$
by choosing a highest weight vector $e_\omega$ and letting inductively  $e_{w-\alpha} = X_{-\alpha}e_w$
for any weight $w$ and any simple root $\alpha$ such that $\omega'-\alpha$ is again a weight. The fact
that $\omega$ is  minuscule implies that the resulting basis is well-defined: $e_w$ does not depend
on the path chosen from $\omega$ to $w$. This follows from \cite[Lemma 1.16]{lw}.
Then
$$y_d(e_w)=\sum_{\langle w,\alpha^\vee\rangle =1}y_d^\alpha e_{w-\alpha}.$$
Since $w-\alpha = s_\alpha(w)$, this should mean on the Schubert
side that
$$p_d=\theta(y_d)=\sum_{\langle w,\alpha^\vee\rangle =1}y_d^\alpha \sigma_{s_\alpha},$$
where $\sigma_{s_\alpha}$ is the Schubert cycle associated to the simple
reflection $s_\alpha$, considered modulo $W_P$.
(This requires a coherent normalization in order to be correct.)
This is remarkably straightforward. In particular we do not need to interpret
$\mathfrak{a}$ in terms of $W$-invariants.

\newtheorem*{one-bis}{Theorem 1 bis}

\begin{one-bis}
There exists an identification of $QH^*(G/P_\omega)$ with $V_\omega^L$, mapping the Schubert
basis $\sigma_\lambda$
to a weight basis $e_\lambda$, such that under this identification,
the quantum multiplication by the special classes $p_d=\theta(y_d)\in H^{2d}(G/P_\omega)$ is given by
$$p_d*\sigma_\lambda=y_{d}(e_\lambda)+qz_{d}(e_\lambda) \;\;\;\; \forall \lambda\in W/W_P.$$

The special classes are identified as follows:

\begin{enumerate}
\item If $G/P_\omega$ is a Grassmannian of type $A_n$, then  $\theta(y_d)=p_d$ is
the class defined by the $d$-th power sum, for $1\le d\le n$.
\item If $G/P_\omega$ is a spinor variety of type $D_n$, then $\theta(y_{2d-1})=p_{2d-1}$
is the class defined by the $(2d-1)$-th power sum, for $1\le d\le n$, and in
degree $d=n-1$ there is also a class $\theta(y'_{n-1})=\tau_{n-1}$.
\item For minuscule $G/P_\omega$ where $G=E_6$ or $G=E_7$, the classes $p_d$ are given explicitely
in section \ref{exceptional} in terms of Schubert classes.
\end{enumerate}
\end{one-bis}

\section{The spectral theorem and the complex involution}

\subsection{Cyclic elements and quantum cohomology}
Recall the Borel presentation of the  cohomology ring of $G/P_{\omega}$,
$$H^*(G/P_\omega,\CC)\simeq \CC[\fh]^{W_P}/(f_{d_1},\cdots ,f_{d_{r-1}},f_{d_r}).$$
Here $d_1,\ldots ,d_r$ are again the exponents of the group $G$
and $f_{d_1},\ldots ,f_{d_{r-1}},f_{d_r}$
are homogeneous generators of $\CC[\fh]^{W}$, of degrees $d_1,\ldots ,d_r$.
The presentation of the quantum cohomology is obtained by deforming the relations
$f_{d_1},\ldots ,f_{d_{r-1}}$, $f_{d_r}$.
When $G/P_\omega$ is minuscule, it turns out that
$$\deg q = d_r=h+1,$$
where $h$ denotes the Coxeter number. In particular only the highest degree
relation $f_{d_r}$ can be deformed, and a suitable normalization yields
$$QH^*(G/P_\omega)\simeq \CC[\fh]^{W_P}[q]/(f_{d_1},\cdots ,f_{d_{r-1}},f_{d_r}-q).$$
(See \cite{cmp1} for more details.)

On the other hand, the cyclic elements in $\fh$
are characterized by the relations $f_{d_1}=\cdots =f_{d_{r-1}}=0, f_{d_r}\ne 0$
\cite{kos}.
Moreover, the value of $f_{d_r}$ completely determines the $W$-orbits of cyclic
elements in $\fh$. Denote by $\mathcal{O}_q$ the $W$-orbit of cyclic
elements on which $f_{d_r}=q$, so that, as observed in \cite{gp}, the ideal
$(f_{d_1},\cdots ,f_{d_{r-1}},f_{d_r}-q)$ is exactly the ideal of $\mathcal{O}_q$.
Kostant has proved that $\mathcal{O}_q$ (for $q\ne 0$) is a free $W$-orbit.
This allows to recover  a result of \cite{cmp3}:

\begin{prop}
If $G/P_\omega$ is minuscule, its quantum cohomology ring
$$QH^*(G/P_\omega)_{q=1}=\CC[\mathcal{O}_1]^{W_P}$$
at $q=1$ is semisimple.
\end{prop}

In other words, $QH^*(G/P_\omega)_{q=1}$ is the algebra of functions on the regular
scheme $Z_{G/P_{\omega}}\simeq \mathcal{O}_1/W_P$ parametrizing $W_P$-orbits in $\mathcal{O}_1$.
We can decompose
$$QH^*(G/P_\omega)_{q=1}=\oplus_{\zeta\in Z_{G/P_{\omega}}}\CC \iota_\zeta,$$
where $\iota_\zeta\in QH^*(G/P_\omega)_{q=1}$ is the primitive idempotent defined
by the point $\zeta$ of $Z_{G/P_{\omega}}$ (the function equal to one at $\zeta$,
and to zero at the other points of  $Z_{G/P_{\omega}}$).

Of course the same would hold for any fixed value of $q$, with a finite
scheme  $Z_{G/P_{\omega}}(q)$ which is nothing else than $Z_{G/P_{\omega}}$, up to
a homothety by some root of $q$.

Note that once we choose an element in $\mathcal{O}_q$, we get an identification
of $Z_{G/P_{\omega}}(q)$ with $W/W_P$. Of course we can choose any Cartan subalgebra in the Borel
presentation. In particular we can choose $\mathfrak{a_q}$. We may suppose that our
prefered cyclic element $x_q$ belongs to $\mathcal{O}_q$. Indeed, $f_{d_r}(x_q)$
depends linearly on $q$, and we may therefore suppose that our root vectors have
been chosen in such a way that $f_{d_r}(x_q)=1$.  This provides us with an
identification  $Z_{G/P_{\omega}}(q)\simeq W/W_P$.

On the other hand, the choice of the highest weight
$\omega$ of $V_\omega$ also identifies the set of weights in $V_\omega$ with $W/W_P$,
which gets therefore identified with $Z_{G/P_{\omega}}(q)$.  This allows to
define
$$\theta_q : \mathfrak{a_q}\rightarrow QH^*(G/P_\omega)_q$$
in the simplest possible way. For any element $y$ of the Cartan algebra
$\mathfrak{a_q}$, $\theta_q(y)$ is the function on $Z_{G/P_{\omega}}(q)$ whose
walue at $\zeta$, considered as a weight of $V_\omega$, is $\zeta(y)$.
In other words,
$$\theta_q(y)=\sum_\zeta \zeta(y)\iota_\zeta.$$
The fact that $\theta_q(x_q)=h$ is the hyperplane class
is purely tautological, once one remembers that, in the Borel presentation,
the hyperplane class is identified with the class of $\omega$, considered as a $W_P$-invariant
function on $\mathfrak{a_q}$.

Moreover the main statement of the theorem is now also essentially tautological.
Indeed, identify $QH^*(G/P_\omega)_q$ with $V_\omega^L$ by sending the idempotent
$\iota_\zeta$ to (any non-zero multiple of) the weight vector $f_\zeta$
of weight $\zeta$ in $V_\omega^L$. Then by construction,
$\theta_q(y)$ acts on $\iota_\zeta$ by multiplication by
$\zeta(y)$, exactly as $y$  acts on the weight vector $f_\zeta$. \qed

\subsection{} However, this does not really clarify the action
of $\mathfrak{a}$ in $QH^*(G/P_\omega)$, which we would like to express in
terms of the Schubert classes $\sigma_\lambda$ rather than on the
idempotent classes $\iota_\zeta$, which don't have much geometrical
meaning. On the other hand, we have
two natural bases in $V_\omega^L$, consisting either of weight
vectors $e_\lambda$ for the original Cartan subalgebra $\mathfrak{h}$,
or of weight vectors $f_\zeta$ for the Cartan subalgebra $\mathfrak{a_q}$.
Of course these two basis are defined only up to scalars.
A remarkable fact, for which we do not have any satisfactory
explanation, is that after suitable normalizations, the vectors
$e_\lambda$ can be expressed in terms of the $f_\zeta$, exactly
as the Schubert classes $\sigma_\lambda$ are expressed in terms
of the idempotents $\iota_\zeta$.

\begin{theo}\label{change}
For any weight $\lambda$ of a classical module $V_\omega^L$,
$$e_\lambda=\sum_\zeta \sigma_\lambda(\zeta)f_\zeta.$$
\end{theo}

We will prove this by an explicit computation in the classical cases, that is
for Grassmannians, even dimensional quadrics, and spinor varieties.
The proof is completely straightforward for Grassmannians and
quadrics, but a bit tricky for spinor varieties, where Schubert
classes are represented by Schur $\tilde{P}$-functions.
These cases being simply-laced we will not really have to deal with Langlands duality.
Note nevertheless that spinor varieties can also be considered as
minuscule spaces of type B, and the latter is exchanged with type C by Langlands duality.
So in principle we should be able to relate the quantum cohomology of
spinor varieties with the biggest fundamental representations of
symplectic groups.

In the two exceptional cases we proceed differentely. First we identify
the special classes $p_d$, using the method we explained above.
Then we use computer programs to compute their action on Schubert classes
by quantum multiplication. Finally we check that this action
coincides with the Lie algebra action of the corresponding elements
of $\mathfrak{e}_6$ and $\mathfrak{e}_7$ on the minuscule representation.
That is, we directly check Theorem 1 bis explicitly,  rather than Theorem \ref{change}.

\subsection{The complex involution.}Another important property is that we could have chosen for cyclic element the
weighted combination
$$w_q=\sum_{i=1}^rn_i(\psi)^{\frac{1}{2}}X_{-\alpha_i}+qX_\psi.$$
This cyclic element has the nice property that its conjugate with respect to
Weyl's compact form is
$$w_q^*=\sum_{i=1}^rn_i(\psi)^{\frac{1}{2}}X_{\alpha_i}+qX_{-\psi},$$
so that $w_q$ is normal in the sense that $[w_q,w_q^*]=0$. In particular the
centralizer of $w_q$ is real, and therefore its graded basis $y_d+qz_d$ is
self-conjugate (up to real scalars). This implies:

\begin{prop}
Let $\iota$ denote the algebra involution of $QH^*(G/P_\omega)_{loc}$ induced
by the complex conjugation. There exists scalars $\gamma_d$ such that
$$\iota(p_d)=q^{-1}\gamma_dp_{h-d}.$$
\end{prop}

Quantum Satake and the complex involution are therefore deeply intertwined.
\medskip

{\bf Acknowledgements.} We thank Michael Finkelberg, Victor Ginzburg,
Yuri Manin and Leonid Rybnikov for discussions
of certain aspects of this paper.

\section{Type A. Grassmannians and power sums}

\subsection{Pl\"ucker coordinates and Schur functions}

The minuscule homogeneous spaces for $G=SL_n$ are the Grassmannians $G(a,b)$, for
$a+b=n$. The quantum cohomology ring of the Grassmann variety $G(a,b)$ is a quotient
of the ring of symmetric functions of $a$ indeterminates, namely
$$QH^*(G(a,b))=\CC [e_1,\ldots , e_a]/(h_{b+1},\ldots , h_{n-1},h_n-q).$$
For symmetric functions our main reference is \cite{mcd}, and we use the same notations.
In particular the $e_r$ and $h_s$ are the elementary and complete symmetric functions,
repectively.

The scheme  $Z_{G(a,b)}$ is the set of (unordered) $a$-tuples $\zeta=(\zeta_1,\ldots ,\zeta_a)$
of $n$-th roots of $q$. The Schubert classes $\sigma_\lambda$ are indexed by partitions
$\lambda =(\lambda_1\ge\cdots \ge\lambda_a)$ such that $\lambda_1\le b$ (we write
$\lambda\subset a\times b$, which means that the diagram of $\lambda$ can be
inscribed inside a rectangle of size $a\times b$). As a function on $Z_{G(a,b)}$,
the class $\sigma_\lambda$ is given by the corresponding Schur function,
$$\sigma_\lambda=\sum_{\zeta\in Z_{G(a,b)}}s_\lambda(\zeta)\iota_\zeta.$$

The fundamental representation associated to $G(a,b)$ is $V_{\omega_a}^L
=\wedge^a\CC^n$. A regular nilpotent $x$ is such that $x(e_i)=e_{i-1}$ with respect
to some basis $e_1,\ldots ,e_n$ of $\CC^n$, with the convention that $e_0=0$.
Moreover, its centralizer $\mathfrak{a}=\langle x,x^2,\ldots ,x^{n-1}\rangle $.
The basis of $\CC^n$ that we have chosen defines a Cartan subalgebra $\fh$ of $\fsl_n$
and induces a root decomposition.  The corresponding cyclic element $x_q$
is then given by the matrix
$$x_q={\small \begin{pmatrix}
0&1 &0&0 & \cdots &0&0&0&0 \\
0&0 &1&0 &\cdots &0&0& 0&0 \\
0&0 &0&1 &\cdots &0&0& 0&0 \\
\cdots   &\cdots &\cdots  &\cdots  &
 \cdots &\cdots &\cdots &\cdots &\cdots  \\
0&0 &0&0 &\cdots &0&1&0 & 0 \\
0&0 &0&0 &\cdots &0&0&1&0 \\
0&0 &0&0 &\cdots &0&0&0&1 \\
q&0 &0&0 &\cdots &0&0&0 &0
\end{pmatrix} }.$$
This cyclic element is regular and semisimple:
the vector
$$f_\zeta=\frac{1}{n}\sum_{k=1}^n\zeta^{k-1}e_k$$
is an eigenvector of $x_q$ for the eigenvalue $\zeta$, when $\zeta^n=q$.
Conversely,
$$e_k=\sum_{\zeta^n=q}\zeta^{-k+1}f_\zeta.$$
Note that $h_n(x_q)=q$, in agreement with the normalization condition $x_q\in\cO_q$.

\medskip
In the wedge power $\wedge^a\CC^n$, we have two natural basis. The first one
is indexed by partitions $\lambda\subset a\times b$: to such a partition, we
associate the decomposable vector
$$e_\lambda=e_{\lambda_1+a}\wedge e_{\lambda_2+a-1}\wedge \cdots \wedge e_{\lambda_a+1}.$$
This is a weight basis for the Cartan algebra $\fh$. On the other hand, we can associate
to a $a$-tuple $\zeta=(\zeta_1,\ldots ,\zeta_a)$ in $Z_{G(a,b)}$ the wedge
product $f_{\zeta_1}\wedge\cdots \wedge f_{\zeta_a}$. In fact this is only
defined up to sign, and we will rather let
$$f_\zeta = \Big(
\prod_{1\le i<j\le a}(\zeta_i-\zeta_j)\Big) f_{\zeta_1}\wedge\cdots \wedge f_{\zeta_a}.$$
This only depends on the unordered $a$-tuple $\zeta$, and gives a weight basis for
the Cartan algebra $\mathfrak{a}_q$.

\begin{prop} For any partition $\lambda\subset a\times b$,
$$e_\lambda=\sum_{\zeta\in Z_{G(a,b)}}s_\lambda(\zeta)f_\zeta.$$
\end{prop}

\proof This is a straightforward computation. Expressing the $e_k$'s in terms
of the $f_\zeta$, we get
$$e_\lambda=\sum_{\zeta\in Z_{G(a,b)}}\det(\zeta_i^{\lambda_j+a-j})_{1\le i,j\le a}
f_{\zeta_1}\wedge\cdots \wedge f_{\zeta_a}.$$
But by the very definition of Schur functions,
$$s_\lambda(\zeta)=\frac{\det(\zeta_i^{\lambda_j+a-j})_{1\le i,j\le a}}{\det(\zeta_i^{a-j})_{1\le i,j\le a}},$$
and the claim follows.\qed

\subsection{Quantum product by power sums}

The centralizer of the cyclic element $x_q$ is simply
$$\mathfrak{a}_{q}=\langle x_{q},x_{q}^2,\ldots ,x_{q}^{n-1}\rangle .$$
The eigenvalues of $x_{q}^\ell$ being the $\ell$-th powers of
the $n$-th roots of $q$, the corresponding function on $Z_{G(a,b)}$
is just
$$\pi(x_{q}^\ell)=\sum_{\zeta\in Z_{G(a,b)}}(\zeta_1^\ell+\cdots +\zeta_a^{\ell})
\iota_\zeta.$$
In other words, $\pi(x_{q}^\ell)$ coincides with the function defined by
the $\ell$-th {\it power sum} $p_\ell$. We deduce the following statement, which
contains an explicit realization of (a special instance of) the Satake isomorphism
in type A:

\begin{coro}\label{coroA}
For $1\le \ell\le n-1$, the element $x_q^\ell$ of $\fsl_n$ acts on the basis
$e_\lambda$ of $\wedge^a\CC^n$, exactly as the power sum class $p_\ell$
acts on the Schubert basis $\sigma_\lambda$ of $QH^*(G(a,b))$.

At the classical level, for $1\le \ell\le n-1$, the element $x^\ell$ of $\fsl_n$
acts on the basis $e_\lambda$ of $\wedge^a\CC^n$, exactly as the power sum
class $p_\ell$ acts on the Schubert basis $\sigma_\lambda$ of $H^*(G(a,b),\CC)$.
\end{coro}

For Schubert classes we can deduce the following statement.

\begin{prop}
Define for a partition $\lambda\subset a\times b$ an element $\sigma_\lambda(x)
\in\mathcal{U}(\fsl_n)$ by the identity
$$\sum_\lambda\sigma_\lambda\sigma_\lambda(x)=
\exp(\sum_{t=1}^{n-1}\frac{p_t x^t}{t!}).$$
Then the action of $\sigma_\lambda(x)$ on $\wedge^a\CC^n$ is given by
the Littlewod-Richardson coefficients:
$$\sigma_\lambda(x)(e_\mu)=\sum_\nu c_{\lambda\mu}^{\nu}e_\nu.$$
\end{prop}

\proof The transition  formulas from power sums to Schubert classes
involve the characters of the symmetric group $\mathcal{S}_n$:
if we denote by $\chi^{\lambda}_{\mu}$ the value of the irreducible
character defined by $\lambda$ on the conjugacy class defined by $\mu$,
we have
$$\sigma_\lambda=\sum_\mu z_\mu^{-1}\chi^{\lambda}_{\mu}p_\mu,
\qquad
p_\mu    =\sum_\lambda \chi^{\lambda}_{\mu}\sigma_\lambda,
$$
where $z_\mu=1^{\alpha_1}2^{\alpha_2}\cdots m^{\alpha_m}
\alpha_1!\alpha_2!\cdots \alpha_m!$ if $r$ appears $\alpha_r$ times
in $\mu$. From the first formula and the previous proposition, we deduce
that the Littlewood-Richardson coefficients are the coefficients of the
action on $\wedge^a\CC^n$ of the following element of the universal
envelopping algebra:
$$\begin{array}{rcl}
\sigma_\lambda(x) & = & \sum_\mu z_\mu^{-1}\chi^{\lambda}_{\mu}
x^{\mu_1}\otimes x^{\mu_2}\otimes\cdots \\
 &=& \sum_\mu \chi^{\lambda}_{\mu} x^{\otimes\alpha_1}\otimes
(\frac{x^2}{2})^{\otimes\alpha_2}\otimes \cdots
\end{array}$$
Multiplying by $\sigma_\lambda$ and summing over $\lambda$, we get
$$\begin{array}{rcl}
\sum_\lambda\sigma_\lambda\sigma_\lambda(x) & = &
\sum_\mu p_\mu x^{\otimes\alpha_1}\otimes
(\frac{x^2}{2})^{\otimes\alpha_2}\otimes \cdots \\
& = &\sum_\alpha (p_1x)^{\otimes\alpha_1}\otimes
(\frac{p_2x^2}{2})^{\otimes\alpha_2}\otimes \cdots \\
& = &\exp(\sum_{t=1}^{n-1}\frac{p_t x^t}{t!}).
\end{array}$$
This concludes the proof. \qed

\medskip\noindent {\it Remark}. This formula could make sense for any
minuscule space $G/P_\omega$ (and also at the quantum level). In general, we have a
homogeneous basis $y_t$ of the centralizer of the regular nilpotent $x$
($t$ being some exponent) and we know how to associate to it a
cohomology class $p_t$. Then the identity (with suitable normalization
for $y_t$), in $H^*(G/P_\omega)\otimes\mathcal{U}(\fg)$,
$$\exp(\sum_t p_ty_t)=\sum_{w\in W_P}\sigma_w\sigma_w(x)$$
defines elements $\sigma_w(x)\in \mathcal{U}(\fg)$, and we could hope
that the coefficients in $\sigma_w(x)(e_v)=\sum_uc_{vw}^u e_u$
are the structure constants, that is,
$$\sigma_w\sigma_v=\sum_uc_{vw}^u \sigma_u.$$

\medskip
We can make Corollary \ref{coroA} completely explicit. If we consider the
action of $x^\ell$ on our basis $e_\lambda$, we get a sum of terms obtained
by changing some $e_{\lambda_i+a-i+1}$ into $e_{\lambda_i+a-i+1+\ell}$. If the resulting wedge product is non
zero, then there is an index $j\le i$ such that $\lambda_{j-1}+a-j+2> \lambda_i+a-i+1+\ell>\lambda_j+a-j+1$.
and we obtain the basis vector $e_\mu$ multiplied by $(-1)^{i-j+1}$, where the partition $\mu$ is such that
$\mu_j=\lambda_i+j-i+\ell$, $\mu_{j+1}=\lambda_{j}+1, \ldots , \mu_{i}=\lambda_{i-1}+1$,
and $\mu_k=\lambda_k$ if $k< j$ or $k >i$. This exactly means that $\mu$ is obtained from
$\lambda$ by adding a border rim of size $\ell$. Morever the height of this border rim
is $j-i-1$.

If we consider the action of $x_{q}^\ell$ on our basis $e_\lambda$, we will get
the terms obtained from the action of $x^\ell$, plus $q$ times the terms obtained
by changing some $e_{\lambda_i+a-i+1}$ into $e_{\lambda_i+a-i+1-n+\ell}$. If the resulting wedge product is non
zero, then there is an index $j\ge i$ such that $\lambda_j+a-j+1> \lambda_i+a-i+1-n+\ell>\lambda_{j+1}+a-j$,
and we obtain the basis vector $e_\nu$ multiplied by $(-1)^{j-i}$, where the partition $\nu$ is such that
$\nu_{i}=\lambda_{i+1}-1, \ldots , \nu_{j-1}=\lambda_{j}-1$, $\nu_j=\lambda_i+j-i-n+\ell$,
and $\nu_k=\lambda_k$ if $k< i$ or $k >j$. This exactly means that $\nu$ is obtained from
$\lambda$ by removing a border rim of size $n-\ell$. Morever the height of this border rim
is $j-i-1$.

We deduce the following quantum multiplication rule by power sums:

\begin{theo}
For any $\ell\le n-1$, the quantum product of a power sum class of degree $\ell$,
by a Schubert class $\sigma_\lambda$ in $G(a,b)$, is given by the formula
$$p_\ell*\sigma_\lambda=\sum (-1)^{h(\mu/\lambda)}\sigma_\mu +(-1)^{a-1}q\sum (-1)^{h(\lambda/\nu)}\sigma_\nu,$$
where partitions $\mu$ in the classical part of the product are deduced from $\lambda$
by adding a border rim of size $\ell$ and height $h(\mu/\lambda)$, while partitions $\nu$
in the quantum correction are deduced from $\lambda$ by suppressing
a border rim of size $n-\ell$ and height $h(\lambda/\nu)$.
\end{theo}

Beware that the height of a border rim is equal to the number of rows it
occupies, {\it minus one}.

\noindent {\it Remark}. This statement can also be proved directly using the quantum
Pieri rules and the fact that
$$p_t=\sum_{r+s=t}(-1)^sre_sh_r.$$
Moreover the classical part (for $q=0$) is in \cite{mcd}, Exercise 11 p. 48.

%

\section{Type D. Quadrics and power sums}


Consider an even dimensional quadric $\QQ^{2n-2}\subset \PP V_{\om_1}$,
where $V_{\om_1}$ is the natural representation of $\fso_{2n}$. Let $\kappa$
denote the quadratic form on $V_{\om_1}$ whose annihilator is $\fso_{2n}$.
We can choose two supplementary isotropic subspaces $E_+$ and $E_-$,
with respective basis $(e_i)_{1\le i\le n}$ and $(e_{-j})_{1\le j\le n}$,
such that $\kappa(e_i,e_{-j})=\delta_{ij}$. In this basis the endomorphisms
in $\fso_{2n}$ are those represented by a matrix which is skew-symmetric
with respect to the second diagonal. Diagonal matrices with this property
form a Cartan subalgebra $\fh$, and lower triangular matrices, a Borel
subalgebra $\fb$. Our cyclic element in $\fso_{2n}$ is then, in the basis
$(e_{-n}, \ldots ,e_{-1}, e_1, \ldots ,e_n)$,
$$x_q={\small \begin{pmatrix}
0&1&\cdots &0&0 & 0&0&\cdots &0&0 \\
0&0&\cdots &0&0 & 0&0&\cdots &0&0 \\
\cdots &\cdots &\cdots &\cdots &\cdots  &
 \cdots &\cdots &\cdots &\cdots &\cdots  \\
0&0&\cdots &0&1 & q&0&\cdots &0&0 \\
0&0&\cdots &0&0 & 0&-q&\cdots &0&0 \\
0&0&\cdots &0&0 & 0&-1&\cdots &0&0 \\
0&0&\cdots &0&0 & 0&0&\cdots &0&0 \\
\cdots &\cdots &\cdots &\cdots &\cdots  &
 \cdots &\cdots &\cdots &\cdots &\cdots  \\
1&0&\cdots &0&0 & 0&0&\cdots &0&-1 \\
0&-1&\cdots &0&0 & 0&0&\cdots &0&0 \end{pmatrix} }.$$

This cyclic element is easy to diagonalize.

\begin{lemm}
For any $\zeta$ such that $\zeta^{2n-2}=(-1)^n4q$,
\begin{eqnarray*}
 \hspace{0mm}f_\zeta &=&
e_{-n}+\zeta e_{-n+1}+\cdots +\zeta^{n-2}e_{-2}+\frac{1}{2}\zeta^{n-1}e_{-1}-\\
 & &
 -2\big(\frac{1}{2}e_{n}-\zeta^{-1} e_{n-1}+\cdots +(-\zeta)^{-n+2}e_{2}+
(-\zeta)^{-n+1}e_{1}\big)
\end{eqnarray*}
is an eigenvector of $x_q$ for the eigenvalue $\zeta$.
Moreover the kernel of $x_q$ is generated by the two isotropic vectors
\begin{eqnarray*}
 \hspace{0mm}f_{0_+} =qe_{-1}-e_1 +q^{\frac{1}{2}}(e_n+e_{-n}), \qquad
 f_{0_-} = qe_{-1}-e_1 -q^{\frac{1}{2}}(e_n+e_{-n}).
\end{eqnarray*}
\end{lemm}

On the other hand, the Chow ring of the quadric $\QQ^{2n-2}$ has rank
$2n$, with a Schubert basis given by the following classes. In degree $k<n-1$,
the class $\sigma_k=h^k$ of a linear section of codimension $k$. In degree
$k=n-1$, the classes $\sigma_{n-1}^+$ and $\sigma_{n-1}^-$ of the two rulings of the
quadric in maximal linear spaces. In degree $k>n-1$, the class $\sigma_k=
\frac{1}{2}h^k$ of a codimension $k$ linear space. The Hasse diagram is
the following:

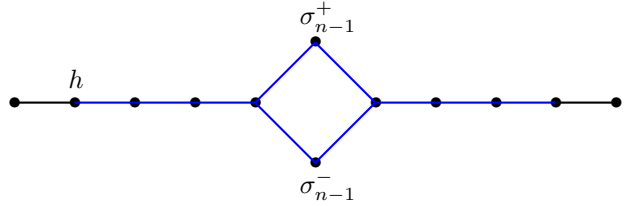
\begin{center}
\psset{unit=4mm}
\psset{xunit=4mm}
\psset{yunit=4mm}
\begin{pspicture*}(20,8)(-10,0)
\multiput(-4.2,3.8)(2,0){5}{$\bullet$}
\multiput(7.8,3.8)(2,0){5}{$\bullet$}
\multiput(5.8,1.8)(0,4){2}{$\bullet$}
\psline(-4,4)(-2,4)
\psline[linecolor=blue](-2,4)(4,4)
\psline[linecolor=blue](4,4)(6,6)
\psline[linecolor=blue](4,4)(6,2)
\psline[linecolor=blue](6,6)(8,4)
\psline[linecolor=blue](6,2)(8,4)
\psline[linecolor=blue](8,4)(14,4)
\psline(14,4)(16,4)
\put(5.5,1){$\sigma_{n-1}^-$}\put(5.5,6.6){$\sigma_{n-1}^+$}
\put(-2.2,4.5){$h$}

\end{pspicture*}
\end{center}

By \cite{cmp3}, the scheme $Z_{\QQ^{2n-2}}$,
for a fixed nonzero value of $q$,
can be identified with the set of complex numbers $\zeta$ such that
$\zeta^{2n-2}=(-1)^n4q$, plus two points that we denote by $0_+$ and $0_-$.
The Schubert classes can be expressed, in terms of the corresponding
idempotents, as
$$\begin{array}{lcl }
\sigma_0 &=& \sum_\zeta \iota_\zeta +\iota_{0_+}+\iota_{0_-},   \\
\sigma_k &=& \sum_\zeta \zeta^k\iota_\zeta, \\
\sigma_{n-1}^+ &=& \frac{1}{2}\sum_\zeta \zeta^{n-1}\iota_\zeta
+q^{\frac{1}{2}}(\iota_{0_+}-\iota_{0_-}),\\
\sigma_{n-1}^- &=& \frac{1}{2}\sum_\zeta \zeta^{n-1}\iota_\zeta
-q^{\frac{1}{2}}(\iota_{0_+}-\iota_{0_-}),\\
\sigma_{n-1+k} &=&  \frac{1}{2}\sum_\zeta \zeta^{n-1+k}\iota_\zeta, \\
\sigma_{2n-2} &=& q\sum_\zeta \iota_\zeta -q\iota_{0_+}-q\iota_{0_-},
\end{array}$$
\noindent where $0<k<n-1$. Quantum Satake then follows
from the observation that the inverse of  the matrix expressing the
$f_\zeta$, $f_{0_+}$ and $f_{0_-}$ in terms of the original basis $e_i$, $e_{-j}$,
is exactly the matrix of the Schubert constants.

\section{Type D. Spinor varieties}

In this section we consider the case of the orthogonal Grassmannian $OG(n,2n)$,
one the two isomorphic connected components of the space of maximal isotropic
subspaces of $\CC^{2n}$, endowed with a non-degenerate quadratic form $\kappa$. Its
minimal equivariant embedding is inside a projectivized half-spin representation
$\PP \Delta$, and for this reason we call it a spinor variety.
The half-spin representation $\Delta$ can be
represented as $\Lambda^{even}E_+$, with the action of $Spin_{2n}$ induced by the
natural action of the Clifford algebra $Cl_{2n}\simeq Cl(E_+)\otimes Cl(E_-)=
\Lambda E_+\otimes \Lambda E_-$. The highest weight line is $1\in \Lambda^0E_+$.
It corresponds to the isotropic space $E_-$, whose family is characterized by the
even dimensionality of the intersection with $E_+$.

\subsection{Weight vectors}
The isotropic spaces corresponding to the weight spaces in $\Delta$ are the
$$E(\epsilon) = \langle   e_{\epsilon_1 1}, \ldots , e_{\epsilon_n n}\rangle,$$
where $\epsilon=(\epsilon_1, \ldots , \epsilon_n)$ is a sequence of signs
such that $\epsilon_1\cdots\epsilon_n=1$. The corresponding line in $\Delta$ is
generated by the vector $v_{E(\epsilon)}$ (a {\it pure spinor}) obtained as the
wedge product of the $e_k$'s  such that $\epsilon_k=1$; there is an even number of such vectors.

If we replace the basis $e_i$, $e_{-j}$ by the basis $f_\zeta$, $f_{0_+}$, $f_{0_-}$
of eigenvectors of $x_q$, the corresponding weight lines in the spin representations
will correspond to maximal isotropic spaces generated by subsets of these basis.
They are obtained by choosing first $n-1$ vectors $f_{\zeta_1}, \ldots , f_{\zeta_{n-1}}$,
where $\zeta_1^2,\ldots,\zeta_{n-1}^2$ are the $n-1$-th distinct roots of $(-1)^n 4q$.
Then we can choose either $f_{0_+}$ or $f_{0_-}$ to generate a maximal isotropic space
-- but note that different choices will produce isotropic spaces meeting in codimension one,
hence belonging to different families (and thus representing lines in different half-spin representations).
We will remain in the same family as $E_-$ if we choose the vector $f_{0_\epsilon}$ and impose the relation
$$\zeta_1 \cdots \zeta_{n-1} = (-1)^n 2\epsilon q^{\frac{1}{2}}$$
for a fixed square root $q^{\frac{1}{2}}$ of $q$. We denote by $F(\zeta)$ the
maximal isotropic space defined by $\zeta_1,\ldots,\zeta_{n-1}$.

What is the corresponding weight line in the half-spin representation?
This has been described in \cite{spinor}. Write the vectors
$f_{\zeta_1}, \ldots ,
f_{\zeta_{n-1}}, f_{0_\epsilon}$ in terms of the original basis. Make a change of
basis giving generators of $F(\zeta)$ of the form $e_-+u(\zeta)e_+$.
The isotropy of $F(\zeta)$ is then equivalent to the fact that the matrix $u(\zeta)$
is skew-symmetric.

\begin{lemm}
The weight line defined by the isotropic space $F(\zeta)$ in the half spin-representation,
  is generated by
$$v_{F(\zeta)}=\sum_\mu \mathrm{Pf}_\mu (u(\zeta))\, e_\mu,$$
where the sum is taken over strict partitions $\mu$ of even length, with
parts taken in $\{1,\ldots ,n\}$.
\end{lemm}

Here $\mathrm{Pf}_\mu (u(\zeta))$ denotes the Pfaffian of the skew-symmetric
matrix obtained by keeping only the rows and columns of $u(\zeta)$ indexed by
the partition $\mu$.

\subsection{Schubert classes}
Schubert classes $\tau_\lambda$ in $OG(n,2n)$ are indexed by strict partitions
$\lambda\subset\rho_{n-1}$,
that is, strictly decreasing sequences $n>\lambda_1>\cdots >\lambda_\ell>0$.
Recall that, in the same way as Schubert classes in ordinary Grassmannians
are related to ordinary Schur functions, Schubert classes in spinor varieties
are related to Schur $\tP$-functions. We briefly recall the definition of these
symmetric functions and the main properties that will be useful to us.
See \cite{pr} for more details.

One starts with a set of variables $x=(x_1,\ldots ,x_{n-1})$, and lets
$$\tP_0(x)=1, \qquad \qquad  \tP_i(x)=e_i(x)/2
\quad \mathrm{for} \;\; 1\le i<n.$$
Then one defines an infinite matrix $\tP(x)=(\tP_{r,s}(x))_{r,s\ge 0}$,
by asking it to be skew-symmetric, and by letting for $r>s$,
$$\tP_{r,s}(x)=\tP_r(x)\tP_s(x)+2\sum_{i=1}^{s-1}(-1)^i\tP_{r+i}(x)\tP_{s-i}(x)+(-1)^{r+s}\tP_{r+s}(x).$$
Finally, let $\lambda$ be any strict partitition. Make it even length if necessary
by adding to it the zero part, and let
\begin{equation}\label{SchurP}
\tP_\lambda(x)=\mathrm{Pf}_\lambda (\tP(x)).
\end{equation}
This is the Schur $\tP$-function associated to the partition $\lambda$. The relation
with Schubert classes in $OG(n,2n)$ is given by the following statement.
Let $\Lambda_{n-1}$ denote the ring of symmetric functions in $n-1$ variables,
with complex coefficients.

\begin{prop}
There is a ring homomorphism $\Lambda_{n-1}\rightarrow H^*(OG(n,2n),\CC)$,
mapping a Schur $\tP$-function $\tP_\lambda$ to the Schubert class
$\tau_\lambda$ if $\lambda$ is a strict partition with parts
smaller than $n$, and to zero otherwise.
\end{prop}

Formula (\ref{SchurP}) then translates into an expression of any Schubert class
in terms of the special ones: this is a {\it Giambelli formula}. It is
proved in \cite{kt} that this formula remains valid in quantum cohomology.

\medskip\noindent {\it Remark}. The Proposition above is the classical version
of Theorem 1 in \cite{kt}. In fact the connection with Schur $P$ and $Q$-functions
has been observed long before. In particular a presentation of $H^*(OG(n,2n),\CC)$
was given by Pragacz in \cite[Theorem 6.17 (ii)]{pr} in terms of $Q$-functions.
Switching from Pragacz's presentation to that of Kresch and Tamvakis
amounts to a formal change of variables.

\medskip
By \cite{cmp3}, the scheme $Z_{OG(n,2n)}$ can be
identified with the set of $(n-1)$-tuples $\zeta=(\zeta_1,\ldots,\zeta_{n-1})$
of complex numbers, such that $\zeta_1^2,\ldots,\zeta_{n-1}^2$ are the
$(n-1)$-th roots of $(-1)^n 4q$. The quantum Giambelli formula then
allows to express any Schubert class as a function on the spectrum,
namely
$$\tau_\lambda=\sum_\zeta \tP_\lambda(\zeta)\,\iota_\zeta.$$

\subsection{The main computation}
Consider a maximal isotropic space $F(\zeta)$ with its basis
$f_{\zeta_1}, \ldots , f_{\zeta_{n-1}}, f_{0_\epsilon}$.
As explained above, in order to compute its representative in
the half-spin representation, we need to express these vectors
as the lines of a matrix $A_+(\zeta)e_+ +A_-(\zeta)e_-$,
and compute the skew-symmetric matrix $u(\zeta)=A_-(\zeta)^{-1}A_+(\zeta)$.

\begin{lemm}
The entries of the skew-symmetrix matrix $u(\zeta)$ are given by
\begin{eqnarray*}
 4 u(\zeta)_{ij} = s_{j-1,1^{i-1}}(\zeta)-s_{j,1^{i-2}}(\zeta).
\end{eqnarray*}
\end{lemm}

\proof The matrices $A_-(\zeta)$ and $A_+(\zeta)$ are
$$A_-(\zeta)=\begin{pmatrix} \frac{1}{2}\zeta_1^{n-1} & \zeta_1^{n-2} & \cdots & \zeta_1 & 1 \\
\cdots & \cdots & \cdots & \cdots & \cdots \\
\frac{1}{2}\zeta_{n-1}^{n-1} & \zeta_{n-1}^{n-2}& \cdots & \zeta_{n-1} & 1\\
q & 0 & \cdots & 0 & \epsilon q^{\frac{1}{2}}
             \end{pmatrix},$$
$$A_+(\zeta)=-2 \begin{pmatrix}  (-\zeta_1)^{-n+1}& (-\zeta_1)^{-n+2} & \cdots & -\zeta_1^{-1} & \frac{1}{2} \\
\cdots & \cdots & \cdots & \cdots & \cdots \\
(-\zeta_{n-1})^{-n+1} & (-\zeta_{n-1})^{-n+2} & \cdots &-\zeta_{n-1}^{-1} & \frac{1}{2}\\
-\frac{\epsilon}{2} q^{\frac{1}{2}} & 0 & \cdots & 0 & \frac{1}{2}.
             \end{pmatrix}$$
First notice that $\det (A_-(\zeta))=(-1)^n2qV(\zeta)$, where $V(\zeta)$ is the
Vandermonde determinant $\det(\zeta_i^{n-j-1})$. Using the usual rules for
the computation of the inverse of a matrix, we deduce that the entries of $u(\zeta)$ can be
obtained as determinants of matrices  equal to $A_-(\zeta)$, up to one column taken
from $A_+(\zeta)$. Namely,
$$ (-1)^{n-1} q V(\zeta)u(\zeta)_{ij}=\det
\begin{pmatrix}  \frac{1}{2}\zeta_1^{n-1} & \zeta_1^{n-2} & \cdots & (-\zeta_1)^{-n+j}& \cdots & \zeta_1 & 1\\
\cdots & \cdots & \cdots & \cdots & \cdots \\
\frac{1}{2}\zeta_{n-1}^{n-1} & \zeta_{n-1}^{n-2} & \cdots & (-\zeta_{n-1})^{-n+j} & \cdots &\zeta_{n-1} & 1\\
q & 0 & \cdots & 0 & \cdots & 0 & \epsilon q^{\frac{1}{2}}
             \end{pmatrix},$$
at least for $1<i,j<n$. Since for all $k$, $(-\zeta_k)^{-n}= \zeta_k^{n-2}/4q$, we can rewrite this, after expanding
with respect to the last row, as
$$\begin{array}{rcl}
4V(\zeta)u(\zeta)_{ij} & = & -\det
\begin{pmatrix} \zeta_1^{n+j-2} & \zeta_1^{n-2} & \cdots & \zeta_1^{n-i+1}& \zeta_1^{n-i-1} &\cdots & \zeta_1 & 1 \\
\cdots & \cdots & \cdots & \cdots & \cdots \\
 \zeta_{n-1}^{n+j-2} & \zeta_{n-1}^{n-2} & \cdots & \zeta_{n-1}^{n-i+1}& \zeta_{n-1}^{n-i-1} &\cdots & \zeta_{n-1} & 1
             \end{pmatrix} \\
 & & +\det \begin{pmatrix} \zeta_1^{n+j-3} & \zeta_1^{n-2} & \cdots & \zeta_1^{n-i}& \zeta_1^{n-i-2} &\cdots & \zeta_1 & 1 \\
\cdots & \cdots & \cdots & \cdots & \cdots \\
\zeta_{n-1}^{2n-j-2} & \zeta_{n-1}^{n-2} & \cdots & \zeta_{n-1}^{n-i}& \zeta_{n-1}^{n-i-2} &\cdots & \zeta_{n-1} & 1
             \end{pmatrix}.
\end{array}$$
Once divided by the Vandermonde, the two determinants of the right hand side are just
Schur functions with hook shapes, as claimed. The computation for $i$ or $j$ equal to
$1$ or $n$ is similar and left to the reader. \qed

\medskip\noindent {\it Remark.} Recall that the elementary symmetric functions $e_r(\zeta)$
and the complete symmetric functions $h_s(\zeta)$ are defined by the formal
identities
\begin{eqnarray*}
e(\zeta;t) &=\sum_{r\ge 0}e_r(\zeta)t^r &= \prod_i(1+t\zeta_i), \\
h(\zeta;t) &=\sum_{s\ge 0}h_s(\zeta)t^s &= \prod_i(1-t\zeta_i)^{-1}.
\end{eqnarray*}
Since $\zeta_1^2,\ldots,\zeta_{n-1}^2$ are the
$(n-1)$-th roots of $(-1)^n 4q$, we have
$$e(\zeta;t)=h(\zeta;t) \times \big(1-(-1)^n 4q t^{2n-2}\big),$$
and in particular $e_r(\zeta)=h_r(\zeta)$ for $r<2n-2$.
This easily implies that
$$s_{k,1^j}(\zeta)=s_{j+1,1^{k-1}}(\zeta)\qquad \mathrm{for}\;j+k<2n-2,$$
and thus that the matrix $u(\zeta)$ is skew-symmetric, as it must be.

\medskip
The hook Schur functions $s_{k,1^j}$ are easy to compute in terms of the
elementary and complete symmetric functions. The Pieri formula
implies the identities
\begin{eqnarray*}
s_{k,1^j} &= &h_{k-1}e_{j+1}-h_{k-2}e_{j+2}+\cdots -(-1)^ke_{j+k}, \\
s_{k,1^j} &= &h_{k}e_{j}-h_{k+1}e_{j-1}+\cdots +(-1)^jh_{j+k}.
\end{eqnarray*}
Using the relation $e_r(\zeta)=h_r(\zeta)$ for $r<2n-2$, we deduce that
for $r+s<2n-2$ and $r>s\ge 0$, we have
$$u(\zeta)_{r+1,s+1}=\frac{e_r(\zeta)}{2}\frac{e_s(\zeta)}{2}+2\sum_{\ell=1}^{s-1}(-1)^\ell
\frac{e_{r+\ell}(\zeta)}{2}\frac{e_{s-\ell}(\zeta)}{2}+(-1)^{r+s}\frac{e_{r+s}(\zeta)}{2}.$$
Formally, this is exactly the expression giving $\tP_{r,s}(\zeta)$.
That means that we can identify $u(\zeta)$ with
a corner of the infinite matrix $\tP(\zeta)$.
We deduce that the coefficients of the vectors $v_{F(\zeta)}$ are
given by values of the Schur $\tP$-functions:
$$v_{F(\zeta)}=\sum_\mu \tP_{\mu-1}(\zeta)e_{\mu}.$$
Note that the parts of $\mu-1=(\mu_1-1,\ldots , \mu_{2l}-1)$ are
strictly smaller than $n$.

Now we can use the orthogonality result for Schur $\tP$-functions
stated in \cite[Proposition 2]{lp}, which implies that for two
weights   $\zeta$ and $\zeta'$ of the half-spin representation,
we have
$$\sum_\lambda \tP_\lambda(\zeta)\tP_{\lambda^c}(\zeta')=
\delta_{\zeta,\zeta'}\prod_i\zeta_i\prod_{j<k}(\zeta_j+\zeta_k),$$
where $\delta_{\zeta,\zeta'}$ is Kronecker's delta.
Therefore, after renormalizing the vectors $v_{F(\zeta)}$ by the nonzero
constants $c(\zeta)=\prod_i\zeta_i\prod_{j<k}(\zeta_j+\zeta_k)$, we can conclude that
the matrix of spectral values of the Schubert classes, is inverse to the
matrix expressing the $v_{F(\zeta)}$'s in terms of the $e_\lambda$'s --
which is precisely what we wanted to prove. \qed

\subsection{The centralizer}
The centralizer of the regular nilpotent element
$x$ in $\fso_{2n}$ is $\mathfrak{a}=\langle x,x^3,\ldots , x^{2n-3}, y\rangle$,
where $y\in\fso_{2n}$ is defined by the conditions that $yx=xy=0$.
Explicitely, we take $y=e_n\otimes e_1^*+e_1\wedge e_n$.
The eigenvalues of $x^{2r-1}$ are the powers of the eigenvalues of $x$,
so that exactly as in type A the corresponding function on
$Z_{OG(n,2n)}$ coincides with the function defined by the power sum $p_{2r-1}$.
On the other hand, the function defined by $y$ associates to a weight
$\eta_1\epsilon_1+\cdots +\eta_n\epsilon_n$, where $\eta_i=\pm 1$,
the sign $\eta_n$. Since there is an even number of minus signs,
$\eta_n=\eta_1\ldots \eta_{n-1}$, and this coincides with
$\zeta_1\ldots \zeta_{n-1}$ up to a constant factor, changing
the signs of the $\eta_k$'s amounts to changing the signs of the
$\zeta_k$'s accordingly.

\subsection{Quantum Satake}
As for ordinary Grassmannians we can translate our results into explicit
product formulas. Beware that we have two different index sets. For
Schubert classes, we used the set of strict partitions $\lambda$
with parts at most $n-1$. Inside $\Delta=\wedge^{even}E_+$, we have
the natural basis $e_\mu$ indexed by strict partitions $\mu$ of even length,
with parts at most $n-1$. We have seen that we must identify these
index sets through the map $\mu\mapsto\lambda=\mu-1$. Conversely,
if $\lambda=(\lambda_1,\ldots ,\lambda_\ell)$ is a strict partition,
with $\lambda_1 <n$, the corresponding weight vector is
$$e^*_\lambda = \Big\{ \begin{array}{ll}
e_{\lambda_\ell+1}\wedge\cdots\wedge e_{\lambda_1+1} & \mathrm{if}\;\ell\; \mathrm{is \; even}, \\
e_1\wedge e_{\lambda_\ell+1}\wedge\cdots\wedge e_{\lambda_1+1} & \mathrm{if}\; \ell\; \mathrm{is\; odd}.
\end{array}.$$
Once this identification has been done, we get the following statement:

\begin{prop}
The product by the class $p_{2r-1}\in QH^*(OG(n,2n))$ defined by a power sum
class of odd degree coincides with the action on the half-spin
representation $\Delta$, of  $x_q^{2r-1}$.
Moreover, the product by the Schubert class $\tau_{n-1}\in QH^*(OG(n,2n))$
coincides with the action on $\Delta$, of
$$y_q=(e_{-n}^*+e_n^*)\otimes (qe_{-1}-e_1)+(e^*_{-1}-qe^*_1)\otimes (e_{-n}+e_n).$$
\end{prop}

\medskip\noindent {\it Remark}.
In terms of Schubert classes, the power sum classes can be expressed as
follows:
$$\begin{array}{rcl}
p_{2r-1} &=& \tau_{2r-1}-2\tau_{2r-2,1}+\cdots +(-1)^{r-1}2\tau_{r,r-1}, \\
p_{2n-2r-1} &=& \tau_{n-1,n-2r}-\tau_{n-2,n-2r+1}+\cdots +(-1)^{r-1}2\tau_{n-r,n-r-1}, \\
\end{array}$$
where $n\ge 2r$ (for $n=2r$ the correct formula for $p_{n-1}$ is the second one). In terms
of the generators $\tau_a$ (corresponding to the elementary symmetric
functions, up to a factor two), since $\tau_a\tau_b=\tau_{a+b}+2\tau_{a+b-1,1}+\cdots
+2\tau_{a-1,b+1}+\tau_{a,b}$ if $a\ge b$,  we have
$$\begin{array}{rcl}
p_{2r-1} &=& (2r-1)\tau_{2r-1}-2(2r-3)\tau_{2r-2}\tau_1+\cdots
+(-1)^{r-1}2\tau_r\tau_{r-1}, \\
p_{2n-2r-1} &=& (2r-1)\tau_{n-1}\tau_{n-2r}-(2r-3)\tau_{n-2}\tau_{n-2r+1}
+\cdots \\
 & & \hspace*{5cm}\cdots+(-1)^{r-1}\tau_{n-r,n-r-1}. \\
\end{array}$$

We can be completely explicit:

\begin{theo}
The quantum product of a Schubert class $\tau_\lambda$
by a power sum class of odd degree is given by
$$p_{2r-1}*\tau_\lambda = \sum (-1)^{h(\mu/\lambda)}\tau_\mu
+2\sum (-1)^{\epsilon(\nu/\lambda)}\tau_\nu
+2q\sum (-1)^{\epsilon(\lambda/\rho)}\tau_\rho.$$
Partitions $\mu$ in the first sum are obtained
by adding to $\lambda$ a border rim of size $2r-1$ and height $h(\mu/\lambda)$.
Partitions  $\nu$ in the second sum are obtained by adding to $\lambda$
a double rim of size $2r-1$, while partitions  $\rho$ in the third sum are
obtained by removing to $\lambda$ a double rim of size $2n-2r-1$.
\end{theo}

The classical part of this formula appears in \cite{mcd}, Exercize 11 p. 265,
to which we refer for the notion of double rims and the definition of
$\epsilon$. We could also obtain the quantum product of any Schubert
class by the special class $\tau_{n-1}$, but one already knows from
\cite[Theorem 6]{bkt} that this product is given by an
extremely simple formula:
$$ \tau_{n-1}\tau_{\lambda}=\tau_{(n-1,\lambda)}+q\tau_{\lambda/(n-1)},$$
where the first (resp. second) term of the right hand side is zero if the
first part of $\lambda$ is equal to (resp. different from) $n-1$.

\section{An example: the Apery equation and $OG(5,10)$}

We illustrate the results of the preceding section by identifying the celebrated differential equation of R.~Apery
with the regularized quantum differential equation for a generic prime Fano threefold of degree twelve.
This must be known to experts but a reference seems to be lacking.

\subsection{Apery's recurrence for  $\zeta(3)$} \rm
Apery proved the irrationality of $\zeta(3)$
in 1979 by considering the recurrence
$$n^3u_n-(34n^3-51n^2+27n-5)u_{n-1}+(n-1)^3 u_{n-2}=0.$$

Denote by $a_n$ he solution of the recurrence with $a_0=1,\; a_1=5$ and
$b_n$, the solution that satisfies $b_0=0,\; b_1=1.$
Then \cite{MP05}

\begin{enumerate}

\item $\displaystyle \left| \zeta(3) - \frac{6 b_n}{a_n} \right| =
\sum_{k=n+1}^\infty \frac{6}{k^3a_ka_{k-1}}= o(a_n^{-2})$;

\item all $a_n$'s are integral; the denominator of $b_n$ divides
$12\mathop{\mathrm{LCM}} (1,2,\dots, n)^3$;

\item $a_n=O(\alpha^n)$ where  $\alpha$ is
 the greatest root of the polynomial
 $x^2-34x+1$.

\end{enumerate}

Put   $\dfrac{6 b_n}{a_n}=\dfrac{p_n}{q_n}$ with
coprime integral $p_n, q_n$. Then it follows from
   $\text{LCM} (1,2, \dots, n)\le (1+\epsilon) e^n$ that
$$\left| \zeta(3) - \frac{p_n}{q_n} \right|= o(q_n^{-1+\delta})$$
for some $\delta > 0$ 
(one can choose $\delta = \frac{\log \alpha -3}{\log \alpha +3}$).
The key assertion here is (2), which follows from the fact
that the solutions  $a_n$ and $b_n$
are iterated binomial sums:
$$a_n= \sum_{k=0}^n \comb{n}{k}^2 \comb{n+k}{k}^2,$$
$$b_n= \frac{1}{6} \sum_{k=0}^n \comb{n}{k}^2 \comb{n+k}{k}^2
\left(\sum_{m=1}^n\frac{1}{m^2}+
\sum_{m=1}^{k} \frac{(-1)^{m-1}}{2m^3\comb{n}{m}\comb{n+m}{m}}\right).$$

\bigskip

We will pass
to the differential operator $L$
(of type D3 in the terminology of \cite{GS07})
that annihilates the generating series $A(t)= \sum a_nt^n, \;B(t)=\sum b_nt^n$
of the Apery numbers. Put
$D=t \dfrac{\partial}{\partial t}$ and
$$L= {D}^{3}-t\left (2\,D + 1\right )\left (17\,{D}^{2}+17\,D+5\right )+{t}^{
2}\left (D+1\right )^{3}
.$$
Then $LA=0$ and $(D-1)LB=0$.

\subsection{Another interpretation of Apery's differential equation: a theorem of Beukers and Peters \cite{BP84}}
  Assume that \mbox{$t \ne 0,1,(\sqrt{2} \pm 1)^4, \infty$}. Then:

\begin{enumerate}
\item The surface $S_t: 1-(1-XY)Z-tXYZ(1-X)(1-Y)(1-Z)=0$
is birationally equivalent to a K3 surface $X_t$;

\item The form
$$\omega_t =  \left. \frac{dX \wedge dZ}{XZ(1-t(1-X)(1-Z)(1-2Y))} \right\vert_{S_t}$$
is the unique (up to scalars) holomorphic 2--form on  $X_t$;

\item $\mathop{ \text{rk Pic}} X_t \ge 19$, with equality for generic $t$;

\item  The periods $y$ of the form $\omega_t$
satisfy the differential equation $Ly=0$.
\end{enumerate}
\rm

\subsection{} \bf A mirror--dual interpretation: varieties $V_{12}$. \rm
A variety $V_{12}$ is a Fano threefold of Picard rank $1$, index $1$ and
anticanonical degree $(-K_V)^3=12$. A generic variety $V_{12}$ can be realized as a section
of the orthogonal Grassmannian $OG(5,10)$
by a linear space of codimension seven \cite{Mukai95}.

We compute its quantum $D$--module and regularized quantum $D$--module
\cite{Golyshev05} via quantum Satake as in the above section,
followed by an iterated application of the quantum Lefschetz principle, see e.g. \cite{CG07} .
Indeed, according to the quantum Chevalley formula, the quantum multiplication matrix
for $OG(5,10)$ is

$$
M =
\left( \begin{array}{cccccccccccccccc}
 {0}&{0}&{0}&{0}&{0}&{0}&{0}&{0}&{0}&{0}&{0}&{t}&{0}&{0}&{0}&{0}\cr
 {1}&{0}&{0}&{0}&{0}&{0}&{0}&{0}&{0}&{0}&{0}&{0}&{0}&{t}&{0}&{0}\cr
 {0}&{1}&{0}&{0}&{0}&{0}&{0}&{0}&{0}&{0}&{0}&{0}&{0}&{0}&{t}&{0}\cr
 {0}&{0}&{1}&{0}&{0}&{0}&{0}&{0}&{0}&{0}&{0}&{0}&{0}&{0}&{0}&{0}\cr
 {0}&{0}&{1}&{0}&{0}&{0}&{0}&{0}&{0}&{0}&{0}&{0}&{0}&{0}&{0}&{t}\cr
 {0}&{0}&{0}&{1}&{0}&{0}&{0}&{0}&{0}&{0}&{0}&{0}&{0}&{0}&{0}&{0}\cr
 {0}&{0}&{0}&{1}&{1}&{0}&{0}&{0}&{0}&{0}&{0}&{0}&{0}&{0}&{0}&{0}\cr
 {0}&{0}&{0}&{0}&{0}&{1}&{1}&{0}&{0}&{0}&{0}&{0}&{0}&{0}&{0}&{0}\cr
 {0}&{0}&{0}&{0}&{0}&{0}&{1}&{0}&{0}&{0}&{0}&{0}&{0}&{0}&{0}&{0}\cr
 {0}&{0}&{0}&{0}&{0}&{0}&{0}&{1}&{1}&{0}&{0}&{0}&{0}&{0}&{0}&{0}\cr
 {0}&{0}&{0}&{0}&{0}&{0}&{0}&{0}&{1}&{0}&{0}&{0}&{0}&{0}&{0}&{0}\cr
 {0}&{0}&{0}&{0}&{0}&{0}&{0}&{0}&{0}&{1}&{0}&{0}&{0}&{0}&{0}&{0}\cr
 {0}&{0}&{0}&{0}&{0}&{0}&{0}&{0}&{0}&{1}&{1}&{0}&{0}&{0}&{0}&{0}\cr
 {0}&{0}&{0}&{0}&{0}&{0}&{0}&{0}&{0}&{0}&{0}&{1}&{1}&{0}&{0}&{0}\cr
 {0}&{0}&{0}&{0}&{0}&{0}&{0}&{0}&{0}&{0}&{0}&{0}&{0}&{1}&{0}&{0}\cr
 {0}&{0}&{0}&{0}&{0}&{0}&{0}&{0}&{0}&{0}&{0}&{0}&{0}&{0}&{1}&{0}\cr
 \end{array} \right)
$$

\medskip
\noindent in the basis of Schubert classes indexed by the  partitions
$\emptyset, (1), (2), (3), (2,1)$, $(4)$, $(3,1)$, $(4,1)$, $(3,2), (4,2), (3,2,1), (4,3), (4,2,1),
(4,3,1), (4,3,2), (4,3,2,1)$.

Choosing $\xi_0=\mathbf{1}$ (the cohomology unit)
for a cyclic vector, we find the minimal quantum differential operator:
$$ [{{\it D}}^{11}\left ({\it D}-1\right )^{5}-t{{\it D}}^{5}\left (2\,
{\it D}+1\right )\left (17\,{{\it D}}^{2}+17\,{\it D}+5\right )+{t}
^{2}] \mathbf{1} =0$$

Now we use quantum Lefschetz to pass to a codimension seven linear section.
It means essentially that, for every $i$
we must multiply the coefficient of $t^i$ (which is in turn a polynomial in $D$) by
$\prod_{j=1}^{i} (D+j)^7$, and strip the result of the trivial parasitic factor on the left.
Thus, we get
$$\begin{array}{l}
    {{\it D}}^{11}\left ({\it D}-1\right )^{5}-t{{\it D}}^{5}\left (2\,
{\it D}+1\right )\left (17\,{{\it D}}^{2}+17\,{\it D}+5\right )(D+1)^7+{t}^{2}(D+1)^7(D+2)^7
  \\
 \hspace*{26mm} =
D^7 (D-1)^5
(D^{4}-t(2D+1)(17D^2+17D+5) +{t}^2(D+1)^2),
 \end{array}$$
so that the quantum differential operator for  $V_{12}$ is
$$D^{4}-t(2D+1)(17D^2+17D+5) +{t}^2(D+1)^2.$$

To arrive finally at the regularized differential operator, we
perform the same multiplication of $t^i$ by $\prod_{j=1}^{i} (D+j)$.
This gives
$$\begin{array}{l}
D^{4}-t(D+1)(2D+1)(17D^2+17D+5) +{t}^2(D+1)^3(D+2) \\
 \hspace*{26mm} = D(D^{3}-t(2D+1)(17D^2+17D+5) +{t}^2(D+1)^3).
\end{array}$$
Stripping away the factor $D$ on the left, we finally get
$$D^{3}-t(2D+1)(17D^2+17D+5) +{t}^2(D+1)^3,$$
which is exactly the operator $L$
above.

\section{Type E. The exceptional minuscule spaces} \label{exceptional}

\subsection{$E_6$ and the Cayley plane}
In order to make explicit computations in $\mathfrak{e}_6$ and its minuscule
representation, we have to choose some basis and be careful with signs.
We start by numbering the simple roots in the following way:

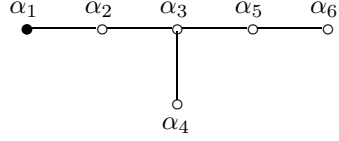
\begin{center}
\setlength{\unitlength}{5mm}
\begin{picture}(20,5)(-5,0)
\multiput(-.3,3.8)(2,0){5}{$\circ$}
\multiput(0,4)(2,0){4}{\line(1,0){1.7}}
\put(3.7,1.8){$\circ$}
\put(-.3,3.8){$\bullet$}
\put(-.6,4.4){$\alpha_1$}
\put(1.4,4.4){$\alpha_2$}
\put(3.4,4.4){$\alpha_3$}
\put(5.4,4.4){$\alpha_5$}
\put(7.4,4.4){$\alpha_6$}
\put(3.45,1.3){$\alpha_4$}
\put(3.85,2.1){\line(0,1){1.8}}
\end{picture}
\end{center}

Then for each simple root $\alpha_i$ we choose a $\fsl_2$-triple
$(X_{-\alpha_i}, H_{\alpha_i}, X_{\alpha_i})$, with $X_{\pm\alpha_i}$
a generator of the corresponding root space in $\mathfrak{e}_6$.
Then the $H_{\alpha_i}=[X_{\alpha_i},X_{-\alpha_i}]$ form a basis of
the Cartan subalgebra. Then we define $X_\beta$ for any root $\beta$,
by induction on the height of $\beta$: if $\mathrm{ht}(\beta)\ge 2$,
let $i$ be maximal such that $\beta-\alpha_i$ is again a root; then
we let
$$X_\beta = [X_{\alpha_i},X_{\beta-\alpha_i}].$$
We can follow the same procedure for negative roots. We end up with a
basis of $\mathfrak{e}_6$, and we can state the following lemma.
We denote by $\theta$ the sum of the simple roots, and by $\psi$
the highest root.

\begin{lemm}
Let $x=X_{\alpha_1}+\cdots +X_{\alpha_6}$. The centralizer $\mathfrak{a}$
of $x$ has a basis given by $y_1=x$, and
$$\begin{array}{rcl}
y_4 & =& X_{\theta-\alpha_2-\alpha_6}+X_{\theta-\alpha_1-\alpha_3}
+X_{\theta-\alpha_1-\alpha_2}-X_{\theta-\alpha_5-\alpha_6}, \\
y_5 & =& X_{\psi-\theta}-2X_{\theta-\alpha_2}
+X_{\theta-\alpha_6}-X_{\theta-\alpha_1},\\
y_7 & =& X_{\theta+\alpha_4}-X_{\theta+\alpha_4+\alpha_3-\alpha_6}+X_{\theta+\alpha_4+\alpha_5-\alpha_1},\\
y_8 & =& X_{\theta+\alpha_4+\alpha_5}-X_{\theta+\alpha_3+\alpha_4}, \\
y_{11} &=& X_{\psi}.
\end{array}$$
\end{lemm}

The basis is indexed by the exponents of $E_6$, $1, 4, 5, 7, 8, 11$.
Each $y_k$ is a linear combination of root vectors associated to roots
of height $k$.

\medskip
Now we turn to the fundamental representation $V_\omega$ of highest weight $\omega=\omega_1$:
this is one of the two minuscule representations, which are dual one to the other.
Its Hasse diagram is as follows, where each dot represents a weight, starting
from $\omega$ on the extreme left. Going to the right through an edge labeled $i$
means that we apply the simple reflexion $s_{\alpha_i}$, or equivalently, that
we substract $\alpha_i$ to the weight.

\begin{center}
\setlength{\unitlength}{3.5mm}
\begin{picture}(40,12)(-8,2)

\put(-6,8){\line(1,0){6}}
\put(0,8){\line(1,1){2}}
\put(0,8){\line(1,-1){4}}

\put(2,10){\line(1,-1){4}}
\put(2,6){\line(1,1){8}}
\put(4,4){\line(1,1){8}}
\put(8,12){\line(1,-1){8}}
\put(10,14){\line(1,-1){8}}
\put(6,10){\line(1,-1){4}}
\put(10,6){\line(1,1){4}}
\put(14,6){\line(1,1){4}}
\put(16,4){\line(1,1){4}}
\put(18,10){\line(1,-1){2}}
\put(20,8){\line(1,0){6}}
\multiput(-6.2,7.7)(2,0){4}{$\bullet$}
\multiput(1.8,9.7)(4,0){5}{$\bullet$}
\multiput(1.8,5.7)(4,0){5}{$\bullet$}
\multiput(3.8,7.7)(4,0){5}{$\bullet$}
\multiput(3.8,3.7)(12,0){2}{$\bullet$}
\multiput(7.8,11.7)(4,0){2}{$\bullet$}
\multiput(21.8,7.7)(2,0){3}{$\bullet$}
\put(9.8,13.7){$\bullet$}
\put(-5.2,7.65){$1$}
\put(-3.2,7.65){$2$}
\put(-1.2,7.65){$3$}
\multiput(.7,6.7)(2,2){2}{$5$}
\multiput(.7,8.7)(2,-2){3}{$4$}

\multiput(2.7,4.6)(2,2){5}{$6$}
\multiput(4.8,8.6)(2,-2){2}{$3$}
\multiput(8.7,6.6)(2,2){3}{$5$}
\multiput(6.7,10.6)(2,-2){3}{$2$}
\multiput(8.7,12.6)(2,-2){5}{$1$}

\multiput(12.7,6.6)(2,2){2}{$3$}
\multiput(14.7,4.6)(2,2){3}{$4$}
\multiput(16.8,8.6)(2,-2){2}{$2$}
\put(20.8,7.65){$3$}\put(22.8,7.65){$5$}\put(24.8,7.65){$6$}
\end{picture}
\end{center}

We claim that we can define a weight basis of $V_\omega$ in the most
simple way. That is, we start with a highest weight vector $e_\omega$,
and then for any weight $\nu$, we define $e_\nu$  by induction on the
height of $\omega-\nu$. To do that, we just need to consider a simple
root $\alpha_i$ such that $\nu+\alpha_i$ is still a weight, and let
$e_\nu = X_{-\alpha_i}e_{\nu+\alpha_i}.$
We claim this is well-defined. Indeed, if there is another simple
root $\alpha_j$ such that $\nu+\alpha_j$ is still a weight, we can
see on the Hasse diagram that $\nu+\alpha_i+\alpha_j$ is also a
weight. Moreover $\alpha_i+\alpha_j$ is never a root, so that
$X_{-\alpha_i}$ and $X_{-\alpha_j}$ commute. Since
$e_{\nu+\alpha_i}=X_{-\alpha_j}e_{\nu+\alpha_i+\alpha_j}$
and $e_{\nu+\alpha_j}=X_{-\alpha_i}e_{\nu+\alpha_i+\alpha_j}$,
we deduce that $X_{-\alpha_i}e_{\nu+\alpha_i}=X_{-\alpha_j}e_{\nu+\alpha_j}$,
hence the claim.

\medskip
Now we use the program outlined in the introduction. We apply the $y_k$'s
to the highest weight vector $e_\omega$ and express the result as a linear
combination in the $e_\nu$. Replacing each $e_\nu$ by the Schubert class
$\sigma_\nu$ we get the following classes:
$$\begin{array}{rcl}
p_1 & = & \sigma_{1}, \\
p_4 & = & \sigma_{5431}-\sigma_{2431}, \\
p_5 & = & \sigma_{52431}-2\sigma_{65431}, \\
p_7 & = & \sigma_{4265431}-\sigma_{3452431}, \\
p_8 & = & \sigma_{54265431}-\sigma_{34265431}, \\
p_{11} & = & \sigma_{24354265431}.
\end{array}$$
Here have used the following notation: each sequence $i_1\ldots i_k$ encodes a reduced
decomposition $s_{\alpha_{i_1}}\cdots s_{\alpha_{i_k}}$ of an element $w$ of the Weyl group,
and $\sigma_{i_1\ldots i_k}$ is the Schubert class associated to the class of $w$
modulo $W_P$. Note that the degrees of the $p_k$'s are respectively $78, 12, 9, 3, 2, 1$.

\medskip
We can state an explicit version of the Satake isomorphism:

\begin{prop}
For any $w\in W/W_P$ and any exponent $k$, we have
$$p_k\cup \sigma_w = y_k(e_w).$$
\end{prop}

\proof For the left hand side  P.E.
Chaput's program \cite{pec} was used.
The right hand side was computed by hand.
Everything fits perfectly. \qed

\medskip\noindent {\it Remark}. Strange duality, defined in \cite{cmp2}, is an
algebra involution on the quantum cohomology algebra of minuscules spaces
(localized at $q$), which sends $q$ to $q^{-1}$ (up to a factor). With the
normalization of \cite[section 4.5]{cmp2}, we get that
$$\iota(p_1)=xq^{-1}p_{11}, \quad \iota(p_4)=q^{-1}p_{8}, \quad
\iota(p_5)=2xq^{-1}p_{7}.$$

\subsection{$E_7$ and the exceptional Freudenthal variety}

We number the simple roots in the following way:

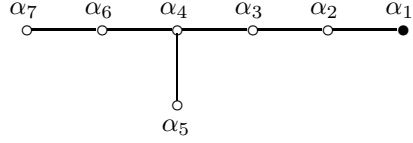
\begin{center}
\setlength{\unitlength}{5mm}
\begin{picture}(20,5)(-5,0)
\multiput(-.3,3.8)(2,0){5}{$\circ$}
\multiput(0,4)(2,0){5}{\line(1,0){1.7}}
\put(3.7,1.8){$\circ$}
\put(9.7,3.8){$\bullet$}
\put(-.6,4.4){$\alpha_7$}
\put(1.4,4.4){$\alpha_6$}
\put(3.4,4.4){$\alpha_4$}
\put(5.4,4.4){$\alpha_3$}
\put(7.4,4.4){$\alpha_2$}
\put(9.4,4.4){$\alpha_1$}
\put(3.45,1.3){$\alpha_5$}
\put(3.85,2.1){\line(0,1){1.8}}
\end{picture}
\end{center}

We deduce a basis of $\mathfrak{e}_7$ by the same procedure we used for $\mathfrak{e}_6$.
The minuscule representation has weight $\omega_1$, its dimension is $56$.
Once we have chosen a highest weight vector, we deduce a basis of $V_{\omega_1}$
by applying the root vectors of the negative simple roots. This produces
the Hasse diagram:

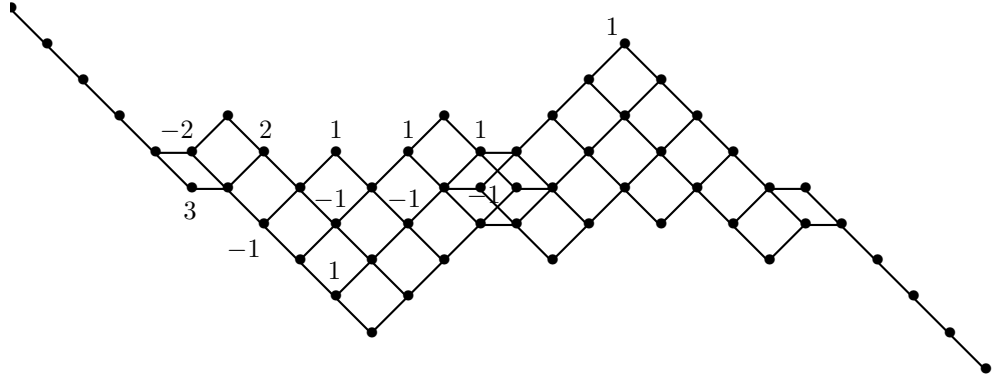
\begin{center}
\psset{unit=2.4mm}
\psset{xunit=2.4mm}
\psset{yunit=2.4mm}
\begin{pspicture*}(50,20)(-8,-6)
\psline(-8,18)(2,8)
\psline(2,8)(4,8)\psline(0,10)(2,10)
\psline(2,10)(4,12)
\psline(2,10)(10,2)\psline(10,2)(12,0)
\psline(4,8)(6,10)
\psline(6,6)(10,10)
\psline(4,12)(12,4)\psline(12,4)(14,2)
\psline(10,10)(14,6)\psline(14,6)(16,4)
\psline(16,12)(18,10)\psline(18,10)(20,8)
\psline(14,10)(16,8)\psline(16,8)(18,6)
\psline(8,4)(16,12)
\psline(10,2)(18,10)
\psline(12,0)(20,8)
\psline(16,8)(18,8)
\psline(18,6)(20,6)
\psline(18,10)(20,10)
\psline(20,8)(22,8)
\psline(18,8)(26,16)
\psline(18,8)(20,6)\psline(20,6)(22,4)
\psline(20,6)(28,14)
\psline(22,4)(30,12)
\psline(20,10)(22,8)\psline(22,8)(24,6)
\psline(22,12)(24,10)\psline(24,10)(28,6)
\psline(24,14)(26,12)\psline(26,12)(34,4)
\psline(26,16)(28,14)\psline(28,14)(36,6)
\psline(28,6)(32,10)
\psline(32,6)(34,8)\psline(34,8)(36,8)
\psline(34,4)(36,6)\psline(36,6)(38,6)
\psline(36,8)(46,-2)
\multiput(-8.3,17.7)(2,-2){6}{$\bullet$}
\multiput(1.7,9.7)(4,0){5}{$\bullet$}
\multiput(5.7,5.7)(4,0){4}{$\bullet$}
\multiput(3.7,7.7)(4,0){5}{$\bullet$}
\multiput(7.7,3.7)(4,0){3}{$\bullet$}
\multiput(3.7,11.7)(12,0){2}{$\bullet$}
\multiput(9.7,1.7)(4,0){2}{$\bullet$}
\put(11.7,-.3){$\bullet$}
\multiput(17.7,7.7)(4,0){5}{$\bullet$}
\multiput(19.7,9.7)(4,0){4}{$\bullet$}
\multiput(19.7,5.7)(4,0){5}{$\bullet$}
\multiput(21.7,11.7)(4,0){3}{$\bullet$}
\multiput(23.7,13.7)(4,0){2}{$\bullet$}
\multiput(25.7,15.7)(4,0){1}{$\bullet$}
\multiput(21.7,3.7)(12,0){2}{$\bullet$}
\multiput(35.7,7.7)(2,-2){6}{$\bullet$}
\put(-9.2,18.3){$1$}\put(25,16.5){$1$}
\put(1.6,6.3){$3$}\put(0.3,10.8){$-2$}
\put(4,4.2){$-1$}\put(5.8,10.8){$2$}
\put(9.7,10.8){$1$}\put(8.8,7){$-1$}\put(9.6,3){$1$}
\put(13.7,10.8){$1$}\put(12.9,7){$-1$}
\put(17.7,10.8){$1$}\put(17.3,7.2){$-1$}
\end{pspicture*}
\end{center}

\medskip
The exponents of $E_7$ are $1, 5, 7, 9, 11, 13, 17$. We denote by $\psi$ the highest root,
and by $\theta$ the sum of the simple roots.

\begin{lemm}
Let $x=X_{\alpha_1}+\cdots +X_{\alpha_7}$. A basis of the centralizer $\mathfrak{a}=\mathfrak{c}(x)$ is
$$\begin{array}{rcl}
y_1 & =& x, \\
y_5 & =& 2X_{\theta-\alpha_1-\alpha_5}+2X_{\theta-\alpha_5-\alpha_7}-X_{\theta-\alpha_1-\alpha_7}+ \\
 & & \qquad
+X_{\theta-\alpha_1-\alpha_2}-X_{\theta-\alpha_1-\alpha_2-\alpha_7+\alpha_4}-3X_{\theta-\alpha_6-\alpha_7}, \\
y_7 & =& 2X_{\theta}-X_{\theta+\alpha_4-\alpha_7}+X_{\theta+\alpha_4-\alpha_1}
-X_{\theta-\alpha_1-\alpha_2+\alpha_4+\alpha_5}+X_{\theta-\alpha_1-\alpha_7+\alpha_3+\alpha_4},\\
y_9 & =& X_{\theta+\alpha_3+\alpha_4}-X_{\theta+\alpha_4+\alpha_5}+X_{\theta+\alpha_2+\alpha_3+\alpha_4-\alpha_7},\\
y_{11} & =& X_{\theta+\alpha_3+2\alpha_4+\alpha_5}-X_{\theta+\alpha_2+\alpha_3+\alpha_4+\alpha_6}
+X_{\theta-\alpha_1+\alpha_3+2\alpha_4+\alpha_5+\alpha_6}, \\
y_{13} &=& X_{\theta+\alpha_2+\alpha_3+2\alpha_4+\alpha_5+\alpha_6}-
X_{\theta+\alpha_2+2\alpha_3+2\alpha_4+\alpha_6}, \\
y_{17} &=& X_\psi.
\end{array}$$
\end{lemm}

We deduce the classes we are interested in:
$$\begin{array}{rcl}
p_1 & = & \sigma_{1}, \\
p_5 & = & 3\sigma_{54321}-2\sigma_{64321}, \\
p_7 & = & 2\sigma_{7564321}-\sigma_{4564321}, \\
p_9 & = & \sigma_{647564321}-\sigma_{347564321}+\sigma_{234564321}, \\
p_{11} & = & \sigma_{43647564321}-\sigma_{23647564321}, \\
p_{13} & = & \sigma_{3243647564321}-\sigma_{2543647564321}, \\
p_{17} & = & \sigma_{432543647564321}. \\
\end{array}$$
Their coefficients on the relevant Schubert classes are indicated on the Hasse diagram above.

\nocite{KKP08}

\end{document}